\documentclass[twoside,english]{amsart}
\advance\oddsidemargin by -1.0cm
\advance\evensidemargin by -1.0cm
\advance\topmargin by -1.0cm 
\usepackage{amssymb}
\usepackage{amsmath}
\usepackage{babel}
\usepackage{amstext}
\usepackage{epsfig}  
\input{diagrams}
\let\mathg\mathfrak
\theoremstyle{plain}
\newtheorem{cor}{Corollary}[section]
\newtheorem{lem}{Lemma}[section]
\newtheorem{thm}{Theorem}[section]
\newtheorem{prop}{Proposition}[section]
\theoremstyle{definition}
\newtheorem{exa}{Example}[section]
\newtheorem{NB}{Remark}[section]
\newtheorem{dfn}{Definition}[section]
\newtheorem*{thank}{Thanks}
\newcommand{\bdm}{\begin{displaymath}}
\newcommand{\edm}{\end{displaymath}}
\newcommand{\ba}[1]{\begin{array}{#1}}
\newcommand{\ea}{\end{array}}
\newcommand{\bea}[1][]{\begin{eqnarray#1}}
\newcommand{\eea}[1][]{\end{eqnarray#1}}
\newcommand{\btab}{\begin{tabular}}
\newcommand{\etab}{\end{tabular}}

\newcommand{\x}{\times}

\newcommand{\ox}{\otimes}

\newcommand{\ra}{\rightarrow}
\newcommand{\lra}{\longrightarrow}

\newcommand{\lmapsto}{\longmapsto}

\newcommand{\tr}{\ensuremath{\mathrm{tr}}}

\newcommand{\del}{\partial}

\newcommand{\ad}{\ensuremath{\mathrm{ad}}}

\newcommand{\C}{\ensuremath{\mathbf{C}}}

\newcommand{\N}{\ensuremath{\mathbf{N}}}
\newcommand{\Z}{\ensuremath{\mathbf{Z}}}

\newcommand{\J}{\ensuremath{\mathcal{J}}}
\renewcommand{\S}{\ensuremath{\mathcal{S}}}

\newcommand{\vphi}{\ensuremath{\varphi}}    

\newcommand{\grad}{\ensuremath{\mathrm{grad\,}}}

\newcommand{\Mor}{\ensuremath{\mathrm{Mor}}}

\newcommand{\Lie}{\ensuremath{\mathrm{Lie}}}

\newcommand{\Ad}{\ensuremath{\mathrm{Ad}\,}}
\newcommand{\Ind}{\ensuremath{\mathrm{Ind}\,}}

\newcommand{\diag}{\ensuremath{\mathrm{diag}}}

\newcommand{\GL}{\ensuremath{\mathrm{GL}}}
\newcommand{\slin}{\ensuremath{\mathg{sl}}}
\newcommand{\SL}{\ensuremath{\mathrm{SL}}}

\newcommand{\so}{\ensuremath{\mathg{so}}}
\newcommand{\SO}{\ensuremath{\mathrm{SO}}}
\newcommand{\Orth}{\ensuremath{\mathrm{O}}}
\newcommand{\p}{\ensuremath{\mathg{p}}}
\renewcommand{\k}{\ensuremath{\mathfrak{k}}}
\newcommand{\g}{\ensuremath{\mathfrak{g}}}
\newcommand{\h}{\ensuremath{\mathfrak{h}}}
\renewcommand{\a}{\ensuremath{\mathfrak{a}}}
\newcommand{\m}{\ensuremath{\mathfrak{m}}}

\newcommand{\X}{\ensuremath{\mathfrak{X}}}
\newcommand{\Y}{\ensuremath{\mathfrak{Y}}}
\newcommand{\nms}{\!\!}
\begin{document}
\setcounter{equation}{0}
%
%
\thispagestyle{empty}
%
\date{\today}
\title[Invariant vector fields on symmetric spaces]{The algebra of 
$K$-invariant vector fields
\\
 on a symmetric space $G/K$}
%
%
%
\author{Ilka Agricola}
\author{Roe Goodman}
\address{\hspace{-5mm} 
{\normalfont\ttfamily agricola@mathematik.hu-berlin.de}\newline
Institut f\"ur Reine Mathematik \newline
Humboldt-Universit\"at zu Berlin\newline
Sitz: WBC Adlershof\newline
D-10099 Berlin\\
Germany}
\address{\hspace{-5mm} 
{\normalfont\ttfamily goodman@math.rutgers.edu}\newline
Rutgers University \newline
Department of Mathematics \newline
110 Frelinghuysen Rd\newline
Piscataway NJ 08854-8019\\
USA}
\thanks{This work was partially supported by the SFB 288 "Differential geometry
and quantum physics" of the Deutsche Forschungsgemeinschaft.}
\begin{abstract}
When $G$ is a complex reductive algebraic group and $G/K$ is a reductive
symmetric space, the decomposition of $\C[G/K]$ as a $K$-module was obtained
(in a non-constructive way) by Richardson, generalizing the celebrated result
of Kostant-Rallis for the linearized problem (the harmonic decomposition of
the isotropy representation). To obtain a constructive version of Richardson's
results, this paper studies the infinite dimensional Lie algebra $\X(G/K)^K$
of $K$-invariant regular algebraic vector fields using the geometry of $G/K$
and the $K$-spherical representations of $G$. Assume $G$ is semisimple and
simply-connected and let $\J$ be the algebra of $K$ biinvariant functions on
$G$. An explicit set of free  generators for the localization
$ \X(G/K)^K_{\psi}$ is constructed for a suitable $\psi \in \J$. A commutator
formula is obtained for $K$-invariant vector fields in terms of the
corresponding $K$-covariant maps from $G$ to the isotropy representation of
$G/K$. Vector fields on $G/K$ whose horizontal lifts to $G$ are tangent to the
Cartan embedding of $G/K$ into $G$ are called \emph{flat}. When $G$ is simple
and simply connected, it is shown that every element of $\X(G/K)^K$ is flat if
and only if $K$ is semisimple. The gradients of the fundamental characters
of $G$ are shown to generate all conjugation-invariant vector fields on $G$.
These results are applied in the case of the adjoint representation of $G =
\SL(2,\C)$ to construct a conjugation invariant differential operator whose
kernel furnishes a harmonic decomposition of $\C[G]$.
\end{abstract}
\maketitle
\pagestyle{headings}
%
%
%
\section{Introduction}\noindent
%
In this paper we study the infinite dimensional Lie algebra of $K$-invariant
vector fields on a reductive symmetric space $G/K$. Our motivation was the
investigation of the algebra of invariant differential operators for non
transitive group actions on smooth affine varieties, and in particular the
abstract Howe duality theorem one has for this situation (see for example
\cite[Satz 2.2]{Agricola01}). Correspondingly, we shall work in the algebraic
category, i.\,e.\ $G$ is a complex connected reductive linear algebraic group
and $K$ is the fixed points of an involutory automorphism $\theta$ of $G$
(thus $G/K$ is the complexification of a Riemannian symmetric space).

There is a canonical $G$-module isomorphism between the space $\X(G/K)$ of
regular algebraic vector fields on $G/K$ and the algebraically induced
representation $\Ind_K^G(\sigma)$, where $\sigma$ is the isotropy
representation of $K$. In particular, the space $\X(G/K)^K$ of $K$-invariant
vector fields on $G/K$ corresponds to the $K$-fixed vectors in the induced
representation. When $G$ is simple and simply connected, Richardson's results
\cite{Richardson82} imply that $\X(G/K)$ is a free module over the algebra
$\J$ of $K$-biinvariant functions on $G$. In Theorem~\ref{fingenvf.thm} we
obtain an explicit set of free generators for a localization
$\X(G/K)^K_{\psi}$, for some $\psi \in \J$.

We next study $\X(G/K)^K$ as a Lie algebra in Section~\ref{liealg.sec}
and obtain a formula for the commutator of $K$-invariant vector fields in
terms of the associated $K$-covariant mappings. The Cartan embedding $G/K
\longrightarrow P \subset G$ given by $gK \mapsto g\theta(g)^{-1}$ is a
fundamental tool in the study of symmetric spaces, and it is natural to use it
to study $\X(G/K)^K$. Invariant vector fields on $G/K$ whose horizontal lifts
to $G$ are tangent to $P$ are called \emph{flat} (in fact, the Cartan
embedding induces {\em a priori} two different notions of flatness, which we
show to be equivalent). We obtain a commutator formula with no curvature term
for the action on $P$ of these vector fields. For $G$ simple and simply
connected, we prove (Theorem~\ref{flatvf.thm}) that every element of
$\X(G/K)^K$ is flat if and only if $K$ is semisimple (i.e. $G/K$ is 
not the complexification of a hermitian symmetric space).

In Section \ref{s-conj} we study the conjugation action of a
semisimple group $G$ on itself. This is an example of the Cartan
embedding of a symmetric space for the group $G\x G$ and involution
$\theta(g,h)=(h,g)$.  In this case all conjugation-invariant vector
fields on $G$ are flat.  Assuming $G$ is simply-connected, we show
that the gradients of the characters of the fundamental
representations of $G$ give a free basis for the conjugation-invariant
vector fields on all of $G$, with no localization needed (Theorem
\ref{conjvf.thm}). In the special case of $\SL(n, \C)$ we calculate
the commutators of an explicit basis of conjugation-invariant vector
fields. When $G=\SL(2,\C)$ we construct a $\C$-basis for $\X_2 =
\X(G)^{\Ad G}$ and compute the commutators and the action on
invariants of this basis. We show that $\X_2$ is isomorphic to a
subalgebra of the Witt algebra (Theorem~\ref{sl2-commutator}) and we
find the highest weight vectors inside $\C[\SL(2,\C)]$.

Section~\ref{s-sep-var-sl2} establishes a separation of variables
theorem for $\SL(2,\C)$. More precisely, using the preceding results,
we construct explicitly a conjugation-invariant differential operator
on $\SL(2,\C)$ such that its kernel $H$ realizes the isomorphism
\bdm
\C[\SL(2,\C)]\ \cong\ \C[\SL(2,\C)]^{\Ad \SL(2,\C)}\ox H \,.
\edm
This result (Theorem~\ref{harmonicity-of-sl2-action}) is the global
version of the separation of variables in the isotropy representation
going back to Kostant and Kostant-Rallis (\cite{Kostant63},
\cite{Kostant&R71}). However, our proof requires extensive
representation-theoretic calculations and does not seem to extend to
arbitrary conjugation actions or symmetric spaces in any obvious way.

A preliminary version of some of the results in this article (in
particular, Theorem~\ref{harmonicity-of-sl2-action}) appeared in the
first author's dissertation (\cite{Agricola00}). While writing this
paper, the first author learned from P.\ Michor (Vienna) that he and
B.\ Kostant had obtained results concerning the conjugation action of
an algebraic reductive group $G$ on itself. In the preprint
\cite{Kostant&M01}, conjugation-equivariant maps and their properties
are studied using an approach somewhat similar to our description via
spherical representations. They also obtain an explicit algebraic
separation of variables theorem for $\SL(n,\C)$.

\begin{thank}
We thank Thomas Friedrich (Humboldt-Universit\"at zu Berlin) for many
valuable discussions on the topic of this paper. We also thank
Siddhartha Sahi (Rutgers University) for his insights and suggestions
concerning spherical representations and Friedrich Knop (Rutgers
University) for helpful comments concerning differential forms and
quotient varieties.
\end{thank}

\section{$K$-Invariant vector fields on $G/K$}\label{s-descr-VF}
\noindent
%
\subsection{Vector fields on $G$ and  $G/K$}
\label{ss-vf}
%
Let $G$ be a connected complex reductive linear algebraic group and
let $\theta$ be an involutive automorphism of $G$. Let $K$ be the
fixed point set $G^{\theta}$. We denote the decomposition of the Lie
algebra $\g$ into the $\pm1$-eigenspaces of $\theta$ by
$\g=\k+\p$. Since $G/K$ is an affine variety, we can identify the
regular functions on $G/K$ with the right $K$-invariant regular
functions on $G$.

We denote by $\X(G)$ (respectively $\X(G/K)$) the regular (algebraic)
vector fields on $G$ (respectively $G/K$). We fix a trivialization of
the tangent bundle $T\, G \cong G \times \g$ so that a regular vector
field $X$ on $G$ corresponds to a regular map $\Phi: G \ra \g$ and a
left-invariant vector field corresponds to a constant map. The
relation between $X$ and $\Phi$ is given by
\begin{equation}
\label{vfield}
   Xf(g) =  \left.\frac{d}{dt}f(g(1+t\Phi(g)))\right|_{t=0}
\end{equation}
for all $f \in \C[G]$ (we may assume $G \subset \GL(n,\C)$; then $f$
is the restriction to $G$ of a regular function on $\GL(n,\C)$ so the
right side makes sense, with the sums and products being matrix
operations). 

\begin{prop}
\label{gkvf.prop}
Let $\sigma$ be the isotropy representation of $K$ on $\p$, and let
$\Ind_K^G(\sigma)$ be the space of regular mappings $\Phi: G \ra \p$
satisfying the right $K$-covariance condition
\begin{equation}
\label{kcovar}
 \Phi(gk)= \sigma(k)^{-1}\Phi(g) \quad \mbox{for all $k\in K$ and $g\in G$.}
\end{equation}
Let $G$ act by left translations on $\Ind_K^G(\sigma)$. Then $\X(G/K)
\cong \Ind_K^G(\sigma)$ as a $G$-module, where the vector field $X \in
\X(G/K)$ corresponding to $\Phi \in \Ind_K^G(\sigma)$ acts by formula
{\rm (\ref{vfield})} on $f \in \C[G/K]$.  In particular, the
$K$-invariant regular vector fields on $G/K$ correspond to the
$K$-fixed elements in $\Ind_K^G(\sigma)$.
\end{prop}

\begin{proof}
The inclusion $\C[G/K] \subset \C[G]$ and the  bundle isomorphism
\bdm
T(G/K)\ \cong\ \ G\x_{K}\p
\edm
imply that a regular vector field $X$ on $G/K$ can be identified with
a map $\Phi \in \Ind_K^G(\sigma)$ by formula (\ref{vfield}). The
covariance condition (\ref{kcovar}) on $\Phi$ implies that $Xf \in
\C[G/K]$ for all $f\in \C[G/K]$; the assumption that $\Phi$ has values
in $\p$ then makes the correspondence $\Phi \lmapsto X$ bijective.
\end{proof}

\subsection{Some finiteness results}
\label{ssection-fin-gen}

Let $\J = \C[ G/K]^K$ be the algebra of $K$-biinvariant regular
functions on $G$. Then $\X(G/K)$ and $\Ind_K^G(\sigma)$ are
$\J$-modules under pointwise multiplication. Furthermore, if $X \in
\X(G/K)$ corresponds to the map $\Phi: G \to \p$, then, for $f \in
\J$, the vector field $fX$ corresponds to the map $g \mapsto
f(g)\Phi(g)$.

Fix a maximal $\theta$-{\em anisotropic} algebraic torus $A\subset G$
with Lie algebra $\a$. Let $M$ be the centralizer of $A$ in $K$, and
$M'$ the normalizer of $A$ in $K$. Let $W = M'/M$ be the ``little Weyl
group''. Then under the restriction map $\J \cong \C[A]^{W}$ (see
\cite[Cor.  11.5]{Richardson82}).

\begin{thm}
\label{invarvf.thm}
Let $G$ be semisimple. Assume that $\C[A]^{W}$ is a polynomial algebra
(this is always true if $G$ is simply connected). Then the space
$\X(G/K)^K$ of $K$-invariant vector fields on $G/K$ is a free
$\J$-module of rank $ \dim (\p^M)$.
\end{thm}

\begin{proof}
We have $\X(G/K)^K \cong \Ind_K^G(\sigma)^K$ as a $\J$-module. But
\bdm
 \Ind_K^G(\sigma)^K = \Mor_K(K\backslash G, \p),
\edm
the space of $K$-equivariant regular maps from the right coset space
$K\backslash G$ to $\p$. By \cite[Theorems 12.3 and
14.3]{Richardson82}, there is a $K$-stable vector subspace $E$ of
$\C[K\backslash G]$ so that the pointwise multiplication map $\J
\otimes E \to \C[K\backslash G]$ is a vector space
isomorphism. Furthermore, the multiplicity of the isotropy
representation $\sigma$ in $E$ is $q = \dim (\p^{M})$. Since $K$ is
reductive, it follows that there are maps $\Phi_1, \ldots, \Phi_q$
from $K\backslash G$ to $\p$ that are linearly independent over $\J$
and span $\Mor_K(K\backslash G)$ as an $\J$-module.
\end{proof}

\begin{NB}
In \cite[Section 15]{Richardson82} Richardson indicates how to
determine all pairs $(G, \theta)$ with $G$ semisimple, such that
$\C[A]^{W}$ is a polynomial algebra.
\end{NB}

\subsection{$K$-invariant vector fields and spherical functions}
\label{ssection-spher-fun}

Given $\psi \in \J$, we use the trivialization of the tangent bundle
of $G$ from Section \ref{ss-vf} to identify the differential of $\psi$
with the map $d\psi: G \to \g^*$ defined by
\bdm
d\psi(g)(X) =  \left.\frac{d}{dt} \psi(g(1+tX))\right|_{t=0}\, 
\quad \mbox{ for $X \in \g$}.
\edm
Since $d\psi(g)(X) = 0$ for $X \in \k$, we can view $g \mapsto
d\psi(g)$ as a map from $G$ to $\p^*$. From the $K$-biinvariance of
$\psi$ it is clear that
\begin{equation}
\label{covar}
  d\psi(kgk')(X) = d\psi(g)(\Ad(k')X)
  \quad  \mbox{for $k, k' \in K$.}
\end{equation}
We fix a bilinear form on $\g$ invariant under $\Ad G$ and
$\theta$. This defines an isomorphism $\p \cong \p^*$ as a $K$-module,
and we let $\grad \psi(g) \in \p$ be the element corresponding to
$d\psi(g) \in \p^*$. From (\ref{covar}) we see that $\grad \psi \in
\Ind_K^G(\sigma)^K$. Hence by Proposition \ref{gkvf.prop} $\grad \psi$
determines a $K$-invariant regular vector field $X_{\psi}$ on $G/K$.

If $(\pi_{\lambda}, V_{\lambda})$ is a finite-dimensional irreducible
$K$-spherical representation of $G$ with highest weight $\lambda$,
then the dual representation $(\pi_{\lambda}^*, V_{\lambda}^*)$ is
also $K$-spherical.  We fix $v_{\lambda}\in V_{\lambda}^K$ and
$v_{\lambda^*} \in V_{\lambda}^{*K}$, normalized so that $\langle
v_{\lambda}, v_{\lambda^*} \rangle = 1$, and let
\bdm 
\psi_{\lambda}(g) = \langle \pi_{\lambda}(g)
v_{\lambda}, v_{\lambda^*}\rangle 
\edm
be the corresponding spherical function on $G$. Then $\psi_{\lambda}
\in \J$ and hence it determines a $K$-invariant regular vector field
that we denote by $X_{\lambda}$.  Recall that when $\g$ is simple and
$(\g, \k)$ is a symmetric pair, then $\k$ is either semisimple or else
has a one-dimensional center (\cite{Helgason78}).

\begin{thm}
\label{fingenvf.thm}
Assume that $G$ is simply connected, $\g$ is simple, and $G/K$ has
rank $r$.  Let $\varphi_1, \ldots, \varphi_r$ be algebraically
independent generators for $\J$, and let $X_1, \ldots, X_r$ be the
corresponding $K$-invariant vector fields on $G/K$. Then there is a
nonzero function $\psi \in \J$ so that the following holds (where
$\J_{\psi}$ and $\X(G/K)^K_{\psi}$ denote localizations at $\psi$):
\begin{enumerate}
\item[(i)] 
If the Lie algebra $\k$ is semisimple, $X_1, \ldots, X_r$ generate the
$\J_{\psi}$-module $ \X(G/K)^K_{\psi}$.
\item[(ii)]
If the center of $\k$ is non-zero and has basis $J$ with $(\ad J)^2 =
- 1$, let $Y_i$ be the vector field corresponding to the map $g
\mapsto (\ad J) \grad\varphi_i(g)$. Then $X_1, \ldots, X_r, Y_1, \ldots, Y_r$
generate  the $\J_{\psi}$-module   $\X(G/K)^K_{\psi}$ .
\end{enumerate}
\end{thm}

\begin{proof}
We use a modification of arguments from \cite[Theorem 8.1]{Steinberg65},
\cite{Sliman} and \cite{Solomon}. We first observe that 
\begin{equation}
\label{flatgrad}
   \grad \varphi(A) \subset \a \quad \mbox{for all $\varphi \in \J$}.
\end{equation}
This is a consequence of the $KAK$ polar coordinate decomposition of
$G$, and holds for any reductive $G$. For the sake of completeness, we
give a proof.  Consider the restricted root space decomposition
\bdm
  \g = \m + \a + \sum_{\alpha} \g_{\alpha} \, ,
\edm
where $\m = \Lie(M)$. We claim that
\begin{equation}
\label{gradorth}
  d\varphi(a)(X) = 0 \quad\mbox{for all $a \in A$ and $X \in \g_{\alpha}$. }
\end{equation}
To prove this, observe that $X + \theta X \in \k$, so
$  d\varphi(a)(X + \theta(X)) = 0$. The left $K$-invariance of
$\varphi$ gives
\begin{eqnarray*}
 0 & = & \frac{d}{dt} \varphi(a + t(X + \theta X )a) |_{t=0}
     = d\varphi(a)(\Ad(a)^{-1}(X+\theta X))
 \\
   & = & d\varphi(a)(a^{-\alpha}X + a^{\alpha}\theta X).
\end{eqnarray*}
Since we already know that
$ d\varphi(a)(a^{\alpha} X + a^{\alpha} \theta X ) = 0$, we conclude that
\bdm
(a^{\alpha} - a^{-\alpha}) d\varphi(a)(X) = 0.
\edm
Thus  (\ref{gradorth}) holds on the dense open set in $A$ where $a^{\alpha}
\neq a^{-\alpha}$, and hence it holds on all of $A$. But  (\ref{gradorth})
implies that $\grad \varphi(A) \subset (\m + \a) \cap \p = \a$, proving
assertion (\ref{flatgrad}).

Now assume $G$ is simply connected. Then the set $\Lambda_{+}$ of
$K$-spherical highest weights is a free monoid generated by dominant weights
$\mu_1, \ldots, \mu_r$ in $\a^{*}$. Let $ \Lambda = \Z \mu_1 + \cdots + \Z
\mu_r $ be the lattice generated by these weights.  For $\lambda \in
\Lambda_{+}$, define the \emph{monomial symmetric function} $m_{\lambda} \in
\C[A]^W$ by
\bdm
   m_{\lambda}(a) = \sum_{\mu \in W\cdot \lambda} a^{\mu} .
\edm
Put a partial order $\prec$ on $\Lambda$ by $\mu \prec \lambda$ if $\lambda -
\mu$ is a sum of positive restricted roots with nonnegative coefficients. If
$\lambda \in \Lambda_{+}$, the spherical function $\psi_{\lambda}$ is
given on $A$ by a character sum of the form
\begin{equation}
\label{sphere-fun}
  \psi_{\lambda}(a) = c_0 m_{\lambda}(a) + 
    \sum_{\substack{\mu \in \Lambda_{+}
     \\ \mu  \prec \lambda}} c_{\mu} m_{\mu}(a) ,
\end{equation}
where $c_0 \neq 0$ (\cite{Vretare76}; see also \cite[Prop. 9.4]{Helgason94}).
When $\lambda = \mu_i$, we write $\psi_{\lambda} = \psi_{i}$.
Let $\omega_i$ be the character $\omega_i(a) = a^{\mu_i}$ of $A$. Then
\bdm
\Omega = \frac{d\omega_1}{\omega_1} \wedge \cdots \wedge 
         \frac{d\omega_r}{\omega_r}
\edm
is a nowhere vanishing top-degree differential form on $A$. By
formula (\ref{flatgrad}) we can write
\begin{equation}
\label{topwedge}
  \left. d\psi_1 \wedge \cdots \wedge d\psi_r \right|_A = f\, \Omega ,
\end{equation}
where $f$ is a regular function on $A$ that we can calculate using the
differentials of $\psi_i|_A$. Set $\rho = \mu_1 + \cdots + \mu_r$. From
formula (\ref{sphere-fun}) we see that
\bdm
  f(a) = ca^{\rho} +  \sum_{\substack{\mu \in \Lambda
     \\ \mu  \prec \rho}} c_{\mu} a^{\mu} ,
\edm
with $c \neq 0$. Hence $f \neq 0$, so we conclude from formulas (\ref{covar})
and (\ref{topwedge}) that $\{ d\psi_1, \ldots, d\psi_r \}$ is linearly
independent on a dense open set in $G$. Now
\bdm 
 d\varphi_1 \wedge \cdots  \wedge d\varphi_r =
 \frac{\partial(\varphi_1,\ldots,\varphi_r)}{\partial(\psi_1,\ldots,\psi_r)} 
   d\psi_1 \wedge \cdots \wedge d\psi_r  \, .
\edm
Since $\varphi_1, \ldots, \varphi_r$ are assumed to generate $\J$, the
Jacobian factor is nonzero. Hence the differentials $d\varphi_1, \ldots,
d\varphi_r$ are also linearly independent on a dense open set in $G$.

When $\k$ is semisimple, $\p$ is an irreducible $K$-module and $\p^M =
\a$ has dimension $r$ by \cite[Prop.~5.14]{Bang-Jensen}. Let
$\{Z_1,\ldots, Z_r\}$ be a set of free generators for $\X(G/K)^K$
given by Theorem \ref{invarvf.thm}. Then there are functions
$\psi_{ij}\in \J$ such that
\[
    X_j = \sum_{i}  \psi_{ij} \, Z_i \, .
\]
Set $\psi = \det[\psi_{ij}]$. Then $\psi \neq 0$, since the vector
fields $X_1, \ldots, X_r$ are linearly independent on a dense open set
of $G/K$. This implies statement (i) of the theorem.

Now assume $\k$ has center $\C J$ with $\ad(J)^2 = - 1$. The vector
fields $Y_i$ in statement (ii) of the theorem are $K$-invariant. Since
$ \p^M = \a \oplus \ad(J)\a$ by \cite[Lemma 5.7 and
Prop.~5.14]{Bang-Jensen}, the vector fields $X_1, \ldots, X_r, Y_1,
\ldots, Y_r$ are linearly independent on a dense open set of $G/K$ by
the argument above. Hence statement (ii) of the theorem follows
from Theorem \ref{invarvf.thm} and the argument used for statement (i).
\end{proof}

\section{Lie algebra structure}\label{liealg.sec}

\subsection{Commutator formula on $G/K$ }
\label{ssection-comm-form}
The symmetric space $G/K$ is the base of a holomorphic principal $K$-fibre
bundle with total space $G$. The canonical connection $Z: TG\ra \k$ on $G$ has
horizontal space
\bdm
T^h_g G\ =\ \{X\in T_g G:\ Z(X)=0 \}\ =\ \{X\in T_g G:\ dL_{g^{-1}}(X)\in\p \}
\ =\ dL_g(\p)
\edm
at the point $g \in G$. Since we are working in the context of linear
algebraic groups, we can take the differential of left and right translation
as usual matrix multiplication; thus we write $dL_g(X) = g\cdot X$ (matrix
product) for $g \in G$ and $X \in \g$.

Let $X,Y$ be vector fields on $G/K$ corresponding to maps $\Phi, \Psi$ in
$\Ind_K^G(\sigma)$, and let $X^*, Y^*$ be their horizontal lifts to
vector fields on $G$. It is clear from the definition of the canonical
connection $Z$ that the horizontal lift $X^*$ of $X$ to a vector field on $G$
is given by formula (\ref{vfield}). If $f$ is any regular function on $G$,
we have by definition 
\bdm
  (X^*Y^*f)(g) = \left. \frac{d^2}{ds\,dt}
      f\big((g+sg\Phi(g))(1+t\Psi(g+sg\Phi(g)))\big)\right|_{s=t=0} \, .
\edm
Taking a first-order Taylor expansion of $\Psi$ to determine the coefficient
of $st$ in the argument of $f$ on the right side of this equation, we find
that
\bdm
  (X^*Y^*f)(g) = \left. \frac{d}{dr}
      f\big(g(1 + rH(g))\big)\right|_{r=0} \, ,
\edm
where $H(g) = \Phi(g)\Psi(g) + d\Psi_g(g\Phi(g))$. Using the same formula
again with the order of $X$ and $Y$ interchanged, we conclude that for any
regular function $f$ on $G$,
\begin{equation}
\label{liftcomm}
   [X^*,Y^*]f(g)=\left. \frac{d}{dr}
     f\big(g(1 + r[\Phi(g), \Psi(g)]
      + r\Phi\circledast\Psi(g))\big)\right|_{r=0} \, ,
\end{equation}
where we have set
\begin{equation}
\label{commform1}
 \Phi \circledast \Psi(g) := d\Psi_g(g\Phi(g)) - d\Phi_g(g\Psi(g)).
\end{equation}
In formula (\ref{liftcomm}) the term $[\Phi(g), \Psi(g)]$ is in $ \k$ and
arises from the curvature of the canonical connection. When $f$ is right
$K$-invariant, however, this term can be omitted and we obtain the commutator
of $X$ and $Y$ as vector fields on $G/K$. Thus we have proved the following.
\begin{prop}
\label{commform1.prop}
Let $X$ and $Y$ be vector fields on $G/K$ corresponding to the maps $\Phi$ and
$\Psi$ in $\Ind_K^G(\sigma)$, respectively. Then the commutator $[X,Y]$
corresponds to the map $\Phi \circledast \Psi$ defined in formula
{\rm (\ref{commform1})}.
\end{prop}

\begin{NB}
Each term on the right side of (\ref{commform1}) satisfies the right
$K$-covariance condition (\ref{kcovar}). Indeed, if $x\in \p$, then 
\bea 
\notag
d\Phi_{gk}(gkx) &=& \left. \frac{d}{dt} \Phi\big(gk(1+tx)\big)\right |_{t=0} =
\left.\frac{d}{dt}\Phi\big(g(1+t\Ad(k)x)k\big)\right |_{t=0}
 \\ 
 \notag 
 &=& \Ad(k^{-1})
d\Phi_g(g\Ad(k)x)
\eea 
by the $K$-covariance property of $\Phi$. Hence 
\bdm
d\Phi_{gk}(gk\Psi(gk)) = \Ad(k^{-1})d\Phi_g(g\Ad(k)\Psi(gk)) =
\Ad(k^{-1})d\Phi_g(g\Psi(g)) 
\edm 
as claimed. Likewise, if $\Phi$ and $\Psi$ are left $K$-invariant, then so is
the map $g \lmapsto d\Phi_g(g\Psi(g))$.
\end{NB}

\begin{NB}
The commutator formula (\ref{commform1}) can also be obtained from Cartan's
structural equation for the canonical connection, using the fact that this
connection is torsion free.
\end{NB}

\subsection{Cartan embedding and flat vector fields }

The Cartan embedding of the symmetric space $G/K$ into $G$ furnishes an
alternate description of vector fields on $G/K$. This will allow us to
discuss the properties of Lie algebra $\X(G/K)^K$ in more detail in
some cases. The algebraic group version of this embedding is treated
in \cite{Richardson82} (see also \cite[Section 11.2.3]{Goodman&W}). We
summarize the results as follows.

\begin{prop}[Cartan embedding]
For $g,y\in G$ the formula $g \star y = g y \theta(g)^{-1}$ defines an action
of $G$ on itself. The orbit of the identity
$ P \ =\ G \star e\ =\ \{g\,\theta(g)^{-1}:\ g\in G\} $
is a closed irreducible subset of $G$ isomorphic to $G/K$ as a $G$-space
(relative to this action).
\end{prop}
\noindent
 This embedding will be denoted by 
$  j:\ G/K\lra P\subset G,\quad gK\lmapsto  g\,\theta(g)^{-1} $.
Thus we have a commutative diagram 
\begin{diagram}[balance,size=2em]
  &  & G/K\\
  & \ruTo & \dTo_{j}\\
 G&\rTo& P\\
\end{diagram}
where the map $G \to P$ is $g \mapsto g\theta(g)^{-1}$ and the map $G
\to G/K$ is $g \mapsto gK$.  The $\star$-action of $K$ on $P$ is the
usual conjugation action. By abuse of notation, we shall often write
$\Ad g$ both for the conjugation action of $G$ on $G$ as well as the
adjoint representation of $G$ on $\g$. We also denote by $\theta$ the
involution on $\g$ as well as on $G$.  At any point $y$ of $P$, one
has the inclusion of tangent spaces $T_y P\subset T_y G$. Set
$\theta_y = (\Ad y)^{-1}\theta$. This is an involution on $\g$, and
we define $\k_y$ and $\p_y$ to be the $\pm 1$-eigenspaces of
$\theta_y$:
 \bdm
 \k_y\ =\ \{X\in\g:\ \theta_y X= +X \},\quad
 \p_y\ =\ \{X\in\g:\ \theta_y X= -X \} \, .
 \edm
Let $\kappa_y$ and $\pi_y$ be the projections on these spaces:
 \bdm
 \kappa_y \ =\ \frac{1}{2}(1+\theta_y),\quad 
 \pi_y \ =\ \frac{1}{2}(1-\theta_y)\,.
 \edm
Then $\p_y$ is exactly the tangent space $T_y P$, viewed as a subspace
of $\g$ via left translation by $y^{-1}$ \cite[Section
11.2.7]{Goodman&W}, and may be realized as
 \bdm
 T_y P \ =\ \p_y\ =\ \{\Ad y^{-1} X - \theta(X):\ X\in\g \}\,.
\edm
The group $K$ permutes the subspaces $\p_y$, leaving
 $\p_e$ is invariant. 
More precisely, $\Ad k$ maps $\p_y$ to $\p_{kyk^{ -1}}$ in an
equivariant way, as follows.
\begin{lem}
\label{kpaction.lem}
The following diagram is commutative:
\begin{diagram}[balance,size=2em]
\p_e & \rTo^{\pi_y} & \p_y\\
\dTo^{\Ad k} && \dTo_{\Ad k}\\
\p_e & \rTo_{\pi_{kyk^{ -1}}} & \p_{\scriptstyle kyk^{-1}}\\
\end{diagram}
\end{lem}
\begin{proof}
If $X$ is in $\p_y$ and $Y = \Ad k X$, then
\begin{eqnarray*}
 \theta(Y) & = &\theta(\Ad k X) = \Ad\theta (k)\, \theta(X)
   =  - \Ad k\, \Ad y X 
\\
 & = & - \Ad (kyk^{-1})\Ad k X  = - \Ad (kyk^{-1})Y.
\end{eqnarray*}
Hence $\Ad k$ maps $\p_y$ to $\p_{kyk^{ -1}}$, as claimed. The
commutativity of the diagram is as easily verified.
\end{proof}

If $\Phi: P \to \g$ is any regular map, then we can define a regular
vector field $\tilde{\Phi}$ on $P$ by
\bdm
  \tilde{\Phi}f(y) = \left. \frac{d}{dt} f(y + ty\pi_y\Phi(y)) \right|_{t=0}
\edm
for $f \in \C[P]$ and $y\in P$. Now assume that $\Phi(y) \in \p$ for
all $y\in P$. Since $\theta(\Phi(y)) = - \Phi(y)$ in this case, we can
write
\bdm
  y \pi_y \Phi(y) = \frac{1}{2} (y\Phi(y) + y\Ad(y^{-1})\Phi(y)) 
     = \{y, \Phi(y)\},
\edm
where $\{a, b\} = (1/2)(ab + ba)$ is the (normalized) anti-commutator
of the matrices $a, b$. This gives the alternate formula
\begin{equation}
\label{anticommvf}
  \tilde{\Phi}f(y) = \left. \frac{d}{dt}f(y + t\{y, \Phi(y)\})\right|_{t=0} 
\end{equation}
for maps $\Phi$ with values in $\p$. If we assume that $\Phi$ is
$K$-equivariant:
\bdm
 \Phi(kyk^{-1}) = \Ad(k) \Phi(y) \quad\mbox{for all $k \in K$ and $y \in P$},
\edm
then a brief calculation (using Lemma \ref{kpaction.lem}) shows that
$\tilde{\Phi}$ is a $K$-invariant vector field on $P$.

\begin{dfn}
The vector field $\tilde{\Phi}$ is said to be  {\em flat}  if
\begin{enumerate}
\item[(i)]
 $\Phi: P \to \p$ is $K$-equivariant
\item[(ii)]
 $\Ad(y) \Phi(y) = \Phi(y) $ for all $y \in P$.
 \end{enumerate}
\end{dfn}
Since $\{y, \Phi(y)\} = y\Phi(y)$ for a flat field, formula (\ref{anticommvf})
becomes
\begin{equation}
\label{flatvf}
  \tilde{\Phi}f(y) = \left. \frac{d}{dt} f(y + ty\Phi(y)) \right|_{t=0}
\end{equation}
in this case.

\begin{lem}
\label{flatvf.lem}
Let $\Phi: P \to \p$ be a regular, $K$-equivariant map.
 The following are equivalent:
\begin{enumerate}
\item[(i)]
 $\tilde{\Phi}$ is flat;
\item[(ii)]
 $\tilde{\Phi}_y \in T_{y}(P)$ for all $y \in P$;
\item[(iii)]
 $\Phi(A) \subset \a$.
\end{enumerate}
\end{lem}   

\begin{proof}
(i) $\Longleftrightarrow$ (ii): 
From the identification of $T_yP$ with a subspace of $\g$, condition (ii) is
equivalent to
\bdm
  \pi_y \Phi(y) = \Phi(y) \quad \mbox{for all $y \in P$.}
\edm
But $\Phi(y) \in \p$ when $y \in P$, so $\theta\Phi(y) = - \Phi(y)$,
and hence 
\bdm
 \pi_y\Phi(y) = \frac{1}{2}(1 + \Ad(y))\Phi(y)\, .
\edm
This gives the equivalence of (i) and (ii).

(i) $\Longrightarrow$ (iii): Let $a\in A$ be a regular element. Then
$\Phi(a) \in \a$ if and only if $\Ad(a)\Phi(a) = \Phi(a)$. In
particular, (i) implies that $\Phi$ maps the regular elements of $A$
into $\a$. Since the regular elements are dense in $A$, this implies
(iii).

(iii) $\Longrightarrow$ (i): Let $a\in A$ and $k \in K$. Set $y =
kak^{-1}$.  Then $\Phi(y) = \Ad(k)\Phi(a)$ by the $K$-covariance
properties of $\Phi$, so $ \Ad(y)\Phi(y) = \Ad(k)\Ad(a)\Phi(a) $. Now
use (iii) and the $K$-covariance again to obtain
\bdm 
 \Ad(y)\Phi(y) = \Ad(k)\Phi(a) = \Phi(y) \, . 
\edm
Since $\Ad(K)A$ is dense in $P$, this equation holds everywhere on
$P$, and hence $\tilde{\Phi}$ is flat.
\end{proof}

\begin{prop}
\label{flatvf.prop}
Let $\X(P)^K_{\rm flat}$ be the set of all flat vector fields on $P$. 

\begin{enumerate}

\item[(i)]
$\X(P)^K_{\rm flat}$ is a $\J$-submodule of $\X(P)^K$.

\item[(ii)]
If $X, Y \in \X(P)^K_{\rm flat}$ correspond to the maps $\Phi,
\Psi$ respectively, then $[X, Y] = Z$, where $Z$ is the flat vector field
corresponding to the map $\Phi \circledast \Psi $. Hence $\X(P)^K_{\rm flat}$
is a Lie subalgebra of $\X(P)^K$.

\end{enumerate}

\end{prop} 

\begin{proof} 
(i): This is obvious from the definition.

(ii): If $\Phi$ and $\Psi$ are any regular maps from $P$ to $\p$, then a
straightforward calculation as in the proof of formula (\ref{liftcomm}) shows
that
$ [\tilde{\Phi}, \tilde{\Psi}] = \widetilde{\Phi\#\Psi} $,
where $\Phi\#\Psi: P \ra \g$ is defined by
\begin{equation}
\label{commform2}
 \Phi\# \Psi(y) = d\Psi_y(\{y, \Phi(y)\}) - d\Phi_y(\{y, \Psi(y)\}) 
+ \frac{1}{2}[\Phi(y), \Psi(y)] \, .
\end{equation}
Note however that $\Phi\# \Psi$ has values in $\g$ rather than $\p$, in
general, so formula (\ref{commform2}) does not define a Lie algebra structure
on the set of regular maps from $P$ to $\p$. The projection onto $T_yP$ of the
$\k$ component $\frac{1}{2}[\Phi(y), \Psi(y)]$ in formula (\ref{commform2}) is
the {\em curvature term}. 

Now assume that $\Phi$ and $\Psi$ correspond to flat vector fields $X$ and
$Y$. Let $a\in A$. Then the pointwise commutator $[\Psi(a), \Phi(a)] = 0$ by
condition (iii) of Lemma \ref{flatvf.lem}. Hence
\bdm
  [\Psi(kak^{-1}), \Phi(kak^{-1}] = \Ad(k)[\Psi(a), \Phi(a)] = 0 
    \quad \mbox{for $k \in K$}
\edm
by $K$-covariance. Since $\Ad(K)A$ is dense in $P$, it follows that $[\Psi(y),
\Phi(y)] = 0$ for all $y\in P$. Hence the curvature term is zero, $\{y,
\Phi(y)\} = y\Phi(y)$, and $\Phi\# \Psi = \Phi \circledast \Psi$. It is clear
that  $\Phi \circledast \Psi$ satisfies condition (iii) of
Lemma \ref{flatvf.lem}, so $[X, Y]$ is a flat vector field.
\end{proof}

\begin{dfn}
Let $X \in \X(G/K)^K$ and let $X^{*}$ be the horizontal lift of $X$ to a
vector field on $G$. The vector field $X$ is said to be {\em horizontally
flat} if $X^{*}_y \in T_y(P)$ for all $y \in P$.
\end{dfn}

If $X$ is a horizontally flat $K$-invariant vector field on $G/K$ and
$f$ is a regular function on $G$ that vanishes on $P$, then $X^{*}f|_P
= 0$ also. Hence $X^{*}$ restricts to a well-defined vector field on
$P$ that we denote by $X^{\natural}$.  If $X$ is defined by a map
$\Phi \in \Ind_K^G(\sigma)^K$, we see from formulas (\ref{vfield}) and
(\ref{flatvf}) that
$   X^{\natural} = \tilde{\Phi} $.
We note that $\Phi$ is uniquely determined by its restriction to $P$,
since $KP$ is dense in $G$, so $X$ is determined by $X^{\natural}$
when $X$ is horizontally flat. Also $\Phi(kyk^{-1}) = \Phi(y)$ for
$k\in K$ and $y \in P$.  Thus, by Lemma \ref{flatvf.lem}, the flatness
of $X$ is equivalent to the condition $\Phi(A) \subset \a$. In this
case, $X^{\natural} \in \X(P)^K_{\rm flat}$.

\begin{prop}
\label{gradflat.prop}
Let $\varphi \in \J$. Then $X_{\varphi}$ is a horizontally flat vector
field.
\end{prop}

\begin{proof}
This follows from formula (\ref{flatgrad}), Lemma \ref{flatvf.lem},
and the remarks just made.
\end{proof}

Let $j^{*}:\C[P] \to \C[G/K]$ be the algebra isomorphism obtained from
the Cartan embedding ($j^{*}f = f\circ j$ for $f \in \C[P]$). Define
the push-forward vector field $j_{*}(X) = {j^{*}}^{-1} \circ X \circ
j^{*}$ for $X \in \X(G/K)$. Then $j_{*}$ gives an isomorphism between
$\X[G/K]^K$ and $\X[P]^K$.  Suppose $X \in \X(G/K)^K$ is defined by a
map $\Phi: G \to \p$. The left $K$-invariance of $\Phi$ and the
isomorphism $G/K \cong P$ given by the Cartan embedding imply the
existence of a regular map $\Psi: G \to \p$ such that
\bdm
 2 \Phi(g) = \Psi(\theta(g)^{-1}g) \quad \mbox{for all $g\in G$.}
\edm
Let $f \in \C[P]$ and $y \in P$. Since $j^{-1}(y^2) = y$ for $y\in P$, we have
\begin{eqnarray}
\label{cartanvf}
   j_{*}(X)f(y^2) &=& X(j^{*}f)(y) = 
   \left. \frac{d}{dt} f\big(y(1 + t \Phi(y))^2y\big) \right|_{t=0} 
   \\
\notag
    &=& \left.  \frac{d}{dt} 
     f\big(y^2(1 + t(\Ad y)^{-1}\Psi(y^2))\big) 
     \right|_{t=0}
\end{eqnarray}
(note that $t \mapsto (1+ t \Phi(y))\theta(1+ t\Phi(y))^{-1}$ is tangent to $t
\mapsto 1 + 2t\Phi(y)$ at $t = 0$). Equation (\ref{cartanvf}) uniquely
determines $j_{*}(X)$, since the map $y \mapsto y^2$ is surjective on $P$.

When $X \in \X(G/K)^K$ is horizontally flat, it determines two vector fields
on $P$, namely $X^{\natural}$ and $j_{*}(X)$. It is evident from equation
(\ref{cartanvf}) that these vector fields are not the same. However, the two
notions of flatness are related as follows.

\begin{lem}
\label{cartanflat.lem}
Let $X \in \X(G/K)^K$. Then $X$ is horizontally flat if and only if $j_{*}(X)
\in \X(P)^K_{\rm flat}$.
\end{lem}

\begin{proof}
Suppose $j_{*}(X) \in \X(P)^K_{\rm flat}$. Then $j_{*}(X) = \tilde{\Psi}$,
where $\Psi: P \to \p$ is a regular $K$-covariant map such that $\Ad(y)\Psi(y)
= \Psi(y)$ for $y\in P$. From equation (\ref{cartanvf}) we see that
\bdm
   2\Phi(y) = \Psi(y^2) \quad \mbox{for all $y \in P$.} 
\edm
It follows that $\Ad(y)\Phi(y) = \Phi(y)$ for $y \in P$, so $X$ is
horizontally flat by Lemma \ref{flatvf.lem}.

Conversely, if $X$ is horizontally flat, then $\Ad(y)\Phi(y) = \Phi(y)$ for
all $y \in P$. Let the map $\Psi$ be as in equation (\ref{cartanvf}). Since
$\Ad(y)\Psi(y^2) = \Psi(y^2)$ for all $y\in P$, we see that $j_{*}(X) =
\tilde{\Psi}$ by equation (\ref{cartanvf}). The right $K$-covariance of $\Phi$
and the surjectivity of the map $y \mapsto y^2$ on $P$ imply that
\bdm 
  \Psi(kyk^{-1}) = \Ad(k) \Psi(y) \quad \mbox{for $k\in K$ and $y \in P$.}
\edm
Thus $\tilde{\Psi} \in \X(P)^K_{\rm flat}$.
\end{proof}

In light of Lemma \ref{cartanflat.lem}, we shall simply use the term {\em flat}
in the rest of the paper to refer either to a horizontally flat vector field
$X \in \X(G/K)^K$ or to an element in $\X(P)^K_{\rm flat}$.

\begin{thm}
\label{flatvf.thm}
Assume that $G$ is simply connected and $\g$ is simple.
The following are equivalent:

\begin{enumerate}
\item[(i)]
$\k  $ is semisimple.
\item[(ii)]
Every $K$-invariant regular vector field on $G/K$ is flat.
\end{enumerate}
Furthermore, when {\rm (ii)} holds, then $\X(P)^K$ = $X(P)^K_{\rm flat}$.

\end{thm}

\begin{proof}
Let $\varphi_1, \ldots, \varphi_r$ be a set of algebraically
independent generators for $\J$. If $\k$ is semisimple, the
$K$-invariant vector fields corresponding to $\grad \varphi_1, \ldots,
\grad \varphi_r$ are a $\J_{\psi}$-module basis for $\X(G/K)^K_{\psi}$
by Theorem \ref{fingenvf.thm}.  These vector fields are flat by
Proposition \ref{gradflat.prop}. Hence all $K$-invariant vector fields
on $G/K$ are flat by Proposition \ref{flatvf.prop} (the property of
flatness is invariant under localization). On the other hand, if $\k$
is not semisimple, then $\ad(J) \grad \varphi_i(A) \not\subset \a$, so
the corresponding vector field $Y_i$ is not flat by Lemma
\ref{flatvf.lem}.  The last statement follows from Lemma
\ref{cartanflat.lem}.
\end{proof}
\begin{NB}
When $\k$ is not semisimple, the space $G/K$ is the complexification of
a hermitian symmetric space. From Theorem \ref{fingenvf.thm} we have a
direct sum decomposition
\bdm
  \X(G/K)^K_{\psi} =  \X +  \Y, 
\edm
where $\X$ is the Lie algebra of flat rational vector fields generated
over $\J_{\psi}$ by the gradient fields $X_i$, while $\Y$ is generated
over $\J_{\psi}$ by the (nonflat) fields $Y_i$. We have not determined
the commutation relations between $X_i$ and $Y_j$.
\end{NB}

We finish this section with an easy example where the $K$-action is
trivial, yielding the Witt algebra of algebraic vector fields on the
one-sphere.
\begin{exa}[Witt algebra]
The (complexified) one-sphere $\C^*$ is, in Cartan's classification, a
symmetric space of type BDI with the following involution,
\bdm
 \C^*\ =\ \SO(2,\C) / S(\Orth(1,\C)\x \Orth(1,\C))\,,
 \quad
 \theta\left[\ba{rr}a & b\\-b&a\ea\right] =
  \left[\ba{rr}a &- b\\b&a\ea\right]\,.
\edm
One checks that $K=\{\mathbf{1},-\mathbf{1}\}$, $P=\SO(2,\C)$ and thus
$\g=\p=\so(2,\C)\cong \C$. Since the $K$-action on $P$ is by
conjugation, it is trivial, so \emph{any} regular map
$P\cong\C^*\ra\p\cong\C$ induces a $K$-invariant vector field on the
sphere. Those are spanned by $f_n(x)=x^n$ for $n\in\Z$, with
differential $(df_n)_x(a)=nax^{n-1}$. Since the projection $\pi_x$ is
trivial, the commutator formula (\ref{commform1}) gives
 \bdm
 f_n \circledast f_m(x) 
   \ =\  
   (df_m)_x(x\cdot f_n(x)) - (df_n)_x(x\cdot f_m(x))
   \ =\
   (m-n)x^{n+m}  
    \ =\ 
  (m-n)f_{m+n}(x)\,.
 \edm
This is the well-known commutator relation of the Witt algebra.
The linear combinations $k_n=f_n-f_{-n}$ and $p_n=p_n+p_{-n}$
satisfy the relations
 \bdm
 k_n \circledast k_m \ =\ (m-n)k_{n+m}-(n+m)k_{m-n}, \quad
 p_n \circledast p_m \ =\ (n-m)k_{n+m}+(n+m)k_{n-m}, 
\edm\bdm
 p_n \circledast k_m \ =\ (n-m)p_{n+m}-(n+m)p_{n-m}\,.
 \edm
The Witt algebra thus carries a $\Z_2$-graduation which we shall encounter
again later (Theorem \ref{sl2-commutator}).
\end{exa}

\section{The conjugation action}\label{s-conj}
%
\subsection{Conjugation-invariant vector fields}

Consider the conjugation action of a connected reductive algebraic group $G$
on itself. It fits into the general scheme by choosing $\bar{G}=G\x G$ with
the involution $\theta(g,h)=(h,g)$. Then $K=\{(g,g):g\in G\}$ is the diagonal
embedding of $G$ in $G\x G$, and the Cartan embedding
 \bdm
j: \bar{G} / K\ =\ (G\x G)/G \lra G\x G, \quad 
 (g,h)K\lmapsto (g,h)\theta(g,h)^{-1}\ =\ (gh^{-1},hg^{-1})
 \edm
realizes $P$ as $\{(g,g^{-1}):g\in G\}$, to which there corresponds
$\p=\{(X,-X): X\in \g\}$ on the Lie algebra side. The regular
functions on $\bar{G}/K$ are of the form $\varphi(g,h) = f(gh^{-1})$,
where $f \in \C[G]$.  In particular, on $P$ the function $\varphi$ is
given by $\varphi(g, g^{-1}) = f(g^2)$ (cf.~the proof of Lemma
\ref{cartanflat.lem}).  

The $K$ action on $P$ is by conjugation in each component, so that we
may restrict attention to the first component. Thus $\C[G] \cong
\C[P]$, where $f \in \C[G]$ gives the function $F(g, g^{-1}) =
f(g)$. Conjugation-invariant algebraic vector fields then correspond
to conjugation-equivariant regular maps from $G$ to $\g$, and we
denote them by $\X(G)^{\Ad G}$.  With this identification, the
spherical functions become the irreducible characters of $G$ and the
representation $\sigma$ becomes the adjoint representation of $G$ on
$\g$. The algebra $\J$ consists of the regular class functions on $G$.

\begin{thm}
\label{conjvf.thm}
Assume $G$ is simply connected and $\g$ is semisimple of rank $r$. Let
$\varphi_1, \ldots, \varphi_r$ be the characters of the fundamental
representations of $G$. Then the vector fields $X_1, \ldots, X_r$ on $G$
corresponding to $\grad \varphi_1, \ldots, \grad \varphi_r$ are a $\J$-module
basis for $\X(G)^{\Ad G}$. Furthermore, all conjugation-invariant vector
fields are flat.
\end{thm}

\begin{proof}
Let $T \subset G$ be a maximal torus. We may take 
\bdm
 A = \{(t,t^{-1}) \,:\, t\in T \}, \quad 
 M = \{(t,t) \,:\, t\in T \}.
\edm
The action of $M$ on $\p$ is equivalent to the adjoint action of $T$
on $\g$, hence $\dim\, \p^M = \dim\, T = r$. By \cite[Theorem
8.1]{Steinberg65} the vector fields $X_1, \ldots, X_r$ are linearly
independent on the set of regular elements of $G$. Hence the function
$\psi$ in Theorem \ref{fingenvf.thm} never vanishes on the set of
regular elements, so its zero set is contained in the set $Q$ of
irregular elements of $G$. But $Q$ is a Zariski closed set of
codimension 3 by \cite[Theorem 1.3]{Steinberg65}. Hence $\psi$ must be
constant.  For $\Phi: G \to \g$ a conjugation equivariant map, we have
$\Ad (y)\Phi(y) =\Phi(y^2y^{-1})=\Phi(y)$ for all $y \in G$. Thus the
vector field $\tilde{\Phi}$ is flat.
\end{proof}

\begin{NB} 
\label{conjvf.rem}
Let $N \cong \C^r$ be the cross-section for the set of regular elements of $G$
constructed in \cite[Theorem 1.4]{Steinberg65}. Then Theorem \ref{conjvf.thm}
applies to any set $\{\varphi_1, \ldots, \varphi_r \}$ of generators for $\J$
if it is known that $\{ d\varphi_1, \ldots, d\varphi_r \}$ is linearly
independent at every point of $N$.
\end{NB}

\subsection{Invariant vector fields on $\SL(n,\C)$}
%

We now apply some of our general results to $\SL(n,\C)$. The same method
applies to other classical groups and symmetric spaces using Theorems
\ref{fingenvf.thm} and \ref{flatvf.thm} and the generators for the invariant
polynomials given in \cite[Section 12.4.2]{Goodman&W}.

\begin{thm}\label{conj-inv-VF-sln}
Let $G = \SL(n, \C)$. Define maps $\Phi_k : G \to \g$ by 
\[
 \Phi_k(g)\ =\ g^k- (1/n) \, \tr(g^k)\cdot\mathbf{1}
 \quad \mbox{for $g \in G$}.
\]
Then $\X(G)^{\Ad G}$ is generated (as a module over $\C[G]^{\Ad G}$) by the
vector fields $\tilde{\Phi}_1, \ldots, \tilde{\Phi}_{n-1}$.
\end{thm}
\begin{proof}
Define $ \varphi_k(g) = (1/k)\, \tr(g^{k})$ for $g \in G$. Then for $X
\in \g$ we calculate that
\bdm
  d\varphi_k(g)(X)  =  \left. \frac{d}{dt} \varphi_k(g(1+tX))\right|_{t=0}
    = \tr(g^k X) = \tr(\Phi_k(g)X).
\edm
Using the trace form to identify $\g$ with $\g^*$, we see that $\grad
\varphi_k = \Phi_k$. The restriction of $\varphi_k$ to the diagonal is
a multiple of the power sum of degree $k$, so $\varphi_1, \ldots,
\varphi_{n-1}$ generate the $G$-invariant regular functions.  The
matrices
\[
 X = \left[\begin{array}{ccccc}
  c_1 & -c_2 & \cdots & (-1)^{n-2}c_{n-1} & (-1)^{n-1} 
  \\
  1   &  0  & \cdots &   0                &  0
  \\
  0   &  1  & \cdots &   0                &  0
  \\
  \vdots & \vdots & \ddots  & \vdots & \vdots   
  \\
  0 &   0   & \cdots &   1               &   0
 \end{array}\right]
\]
give a cross-section $N$ for the regular elements of $G$ as $[c_1, c_2,
\ldots, c_{n-1}]$ ranges over $\C^{n-1}$ \cite[Section 7.4]{Steinberg65}. It
is easy to see that $X, X^2, \ldots, X^{n-1}$ are linearly independent. Hence
the maps $\Phi_{1}, \ldots, \Phi_{n-1}$ are linearly independent at all points
of $N$. The result now follows from Remark~\ref{conjvf.rem}.
\end{proof}

\noindent
We compute the commutators of the vector fields in Theorem
\ref{conj-inv-VF-sln}. Since all the conjugation invariant vector
fields are flat (by Theorem \ref{conjvf.thm}), it suffices by Proposition
\ref{flatvf.prop} to calculate the maps $\Phi_k \circledast \Phi_l$.

\begin{thm}
The maps $\Phi_k$ satisfy the 
commutation relations
\bdm
\Phi_k \circledast \Phi_l(g) \ =\ (l-k)\cdot\Phi_{k+l}(g) 
 +\frac{k}{n}\cdot \tr(g^l)\,\Phi_{k}(g) 
- \frac{l}{n}\cdot\tr(g^k)\,\Phi_l(g) \,.
\edm
\end{thm}
\begin{proof}
One obtains for the differential
\bdm
(d\Phi_k)_{g}(X)\ =\ Xg^{k-1}+gXg^{k-2}+ \cdots + g^{k-1}X -
\frac{k}{n}\tr(Xg^{k-1})\,,
\edm
which implies
\bea[*]
 (d\Phi_k)_{g}\big(g\cdot \Phi_l(g)\big)& = & 
 k\bigg(g^{l+1}-\frac{1}{n}\, \tr(g^l)\,g\,\bigg)\,  g^{k-1}
  -\frac{k}{n}\,\tr\bigg(\big(g^{l+1}-\frac{1}{n}
\tr(g^l)g\big)g^{k-1}\bigg)
 \\
 & = & k\big(g^{k+l}-\frac{1}{n}\,\tr(g^{k+l})\, \big)
  - \frac{k}{n}\,\tr(g^l)\big(g^k-\frac{1}{n}\,\tr(g^k) \big)
 \\
 & = & k\cdot \Phi_{k+l}-  \frac{k}{n}\,\tr(g^l)\Phi_k \, .
\eea[*]
Now apply formula (\ref{commform1}).
\end{proof}
\noindent
In particular, the relation
\bdm
\Phi_1 \circledast \Phi_{-1}(g) \ =\ 
 \frac{1}{n} \big(\tr(g)\,\Phi_{-1}(g) +\tr(g^{-1})\,\Phi_{1}(g) \big)
\edm
shows that $\Phi_1$ and $\Phi_{-1}$ generate a finite Lie ring over the ring
of invariants.

\subsection{Invariant vector fields on $\SL(2,\C)$}
%
We  consider the case  $G=\SL(2,\C)$ in more detail.

\begin{thm}\label{conj-inv-VF-sl2}
Every conjugation invariant map $\Psi:\ G=\SL(2,\C) \ra \g=\slin(2,\C)$
is a multiple of the map $\Psi_1: g\mapsto g-g^{-1}$ by an element of
$\C[G]^{\Ad G}$.
\end{thm}
\begin{proof}
The representation of $G$ on $\C^2$ is self-dual, so its character $\chi$
satisfies
\bdm
2 \chi(g) = \tr(g + g^{-1})
\edm
Hence $2 d\chi = \Psi_1$ by the calculation in the proof of Theorem
\ref{conj-inv-VF-sln}. The result now follows from Theorem \ref{conjvf.thm}.
\end{proof}
\noindent
In order to get a $\C$-basis of the space $\X\big(\SL(2,\C)\big)^{\Ad
\SL(2,\C)} = \X_2$, it thus suffices to choose any convenient basis of the
space of invariants. The traces on symmetric tensor powers of the fundamental
representation $V$ of $G$ turned out to yield the simplest formulas.
\begin{prop}\label{sl2-invs-maps}
Let $g$ be an element of $G=\SL(2,\C)$ and denote by $S^k V$ the
$(k+1)$-dimensional irreducible representation of $G$. Then
\bdm
g^{k+1} - g^{-k-1}\, =\, \tr(g)\big|_{S^kV}\cdot (g-g^{-1}) .
\edm
Furthermore,
\bdm
\tr(g)\big|_{S^k V}\, =\, \tr(g^{k})+\tr(g^{k-2})+ \cdots +
\left\{\ba{ll}1 & k\, \text{even}\\ \tr(g) & k\, \text{odd}. \ea\right.
\edm
\end{prop}
\begin{proof}
We first prove the second formula on the maximal torus $T$ of $G$,
chosen as before.  For $h=\diag(x,1/x)\in T$, one has
 \bdm
\tr(h)\big|_{S^k V}\, =\, x^k+x^{k-2} + \cdots + x^{2-k}+ x^{-k}\,.
 \edm
Since $\tr(h^n)=x^n+x^{-n}$, the formula follows immediately on $T$.
The case distinction for the last term depends on whether the number
of summands is even or odd. Since the trace is conjugation invariant,
the formula is valid on the dense set of all conjugates of $T$,
and therefore also holds on $G$. For the first formula, we note
that 
\bdm
\left[\ba{cc}a & b\\ c & d\ea\right] +\left[\ba{cc}d & -b\\ -c & a\ea\right] 
\ =\ (a+d)\cdot \mathbf{1}
\edm
implies $g^n+g^{-n}=\tr(g^n)\cdot\mathbf{1}$, so the algebraic
identity
\bdm
g^{k+1}-g^{-k-1}\ =\ (g^k+g^{k-2} + \cdots + g^{2-k}+g^{-k})\cdot (g-g^{-1})
\edm
finishes the proof.
\end{proof}
\begin{cor}
The vector fields defined by the maps $\Psi_k(g) = g^k - g^{-k}$ for
$k\geq 1$ are a basis for $\X_2$ as a vector space over $\C$.
\end{cor}
\noindent
We compute the commutation relations for this basis.  For notational
simplicity, we write $\Psi_k$ for the conjugation invariant vector
field defined by the map $g  \mapsto \Psi_k(g)$.
\begin{thm}\label{sl2-commutator}
The vector fields $\Psi_k$ satisfy the commutator relations
\bdm
[\Psi_k,\Psi_l]\ =\ (l-k)\Psi_{k+l} - (k+l)\Psi_{l-k}\,.
\edm
In particular, the algebra $\X_2$ of conjugation invariant vector fields
on $\SL(2,\C)$ is isomorphic to a subalgebra of the Witt algebra,
and the vector fields with even index $\{\Psi_{2k}\}_{k\geq 1}$ span
a subalgebra of $\X_2$.
\end{thm}
\begin{proof}
We compute the differential
\begin{eqnarray*}
(d\Psi_k)_{g}(X) &=&  Xg^{k-1} + gXg^{k-2} + \cdots + g^{k-1}X 
\\
  & & \quad  +\ g^{-1}Xg^{-k} + g^{-2}Xg^{-k+1} + \cdots + g^{-k}Xg^{-1}\,,
\end{eqnarray*}
from which we obtain
\bdm
(d\Psi_k)_{g}\big(g\cdot \Psi_l(g)\big)\ =\ k\cdot (g^{k+l}-g^{-k-l}+
g^{l-k}-g^{k-l})\ =\ k\cdot(\Psi_{k+l}(g)+\Psi_{l-k}(g))\,.
\edm
The commutator formula (\ref{commform1}) and Proposition
\ref{flatvf.prop} then imply the result.
\end{proof}
\noindent
The action of the vector fields $\Psi_k$ on the invariants is of
particular interest. The three most important  bases for the invariant
functions are:
 \bdm
 \chi_m(g)\ =\ \tr(g)\big|_{S^m V},\quad I_m(g)\ =\ \tr(g^m),\quad
 J_m(g)\ =\ \tr(g)^m .
 \edm
Only the action of $\Psi_k$ on the power sum $I_m$ is given by a
simple formula.  For this reason, we restrict our attention to $k=1$
in the other two cases.
\begin{thm}\label{action-on-invariants}
The vector field $\Psi_k$ acts on  invariants in 
$\C[\SL(2,\C)]^{\Ad \SL(2,\C)}$ as follows:
\bdm
\Psi_k(I_m)\ =\ m(I_{m+k}-I_{m-k}),
\edm
\bdm
\Psi_1\,(\chi_m)\ =\ m\,\chi_{m+1}-(m+2)\,\chi_{m-1},\quad
\Psi_1(J_m)\ =\ m\, (J_{m+1}-4\, J_{m-1})\,.
\edm
\end{thm}
\begin{proof}
For the invariant $I_m$, the computation is straightforward using
formula~(\ref{vfield}) 
\bea[*]
 \Psi_k(I_m) & =& 
 \left.  \frac{d}{dt}I_m(g + tg\Psi_k(g))\right|_{t=0}
\, =\, \left.\frac{d}{dt}\tr\big[(g+t(g^{1+k}-g^{1-k}))^m \big]\right|_{t=0}
\\
 & = & \,
 \left. \frac{d}{dt}\tr\big[g^m+t\,m\,g^{m-1}(g^{1+k}-g^{1-k})
    + \cdots \big]\right|_{t=0}
\\ 
 & = & \,
 \left. \frac{d}{dt}\tr\big[g^m+t\,m(g^{m+k}-g^{m-k})\big]\right|_{t=0} \\
 & = & 
 m\left(\tr(g^{m+k}) - \tr(g^{m-k})\right)
  \ =\ m(I_{m+k}-I_{m-k})\,.
\eea[*]
Since $\chi_m=I_m+I_{m-2} + \cdots$ by Proposition \ref{sl2-invs-maps}, the
second formula is easily proved by induction. The last formula  is shown
with a similar argument than the first and requires at one stage the identity
$\tr(g^2)=\tr(g)^2-2$, which is immediately verified on matrices.
\end{proof}
\noindent
From the point of representation theory, the infinite dimensional Lie
algebra $\X_2$ comes with two natural representations (and in fact many
more, see Section~\ref{s-sep-var-sl2}). The commutator formula
(Theorem~\ref{sl2-commutator}) describes the structure of the adjoint
representation of $\X_2$ and shows in particular that it has no (non-trivial)
finite-dimensional subalgebras. The action on invariants contains a trivial
summand (the constant function, annihilated by all $\Psi_k$), the rest is
indecomposable in the following sense: for any fixed $m\neq 0$, the linear
hull $V_1$ of the invariants
 \bdm
 I_m,\, I_{m\pm 2}, I_{m\pm 4},\,\ldots, I_{m+1}-I_{m-1},\, 
 I_{m\pm 3}-I_{m\pm 1},\, I_{m\pm 5}-I_{m\pm 3},\,\ldots
 \edm
is invariant under the action of $\X_2$, but its complement $V_2$ spanned by
 \bdm
 I_{m+1}+I_{m-1},\, I_{m\pm 3}+I_{m\pm 1},\, I_{m\pm 5}+I_{m\pm 3},\,\ldots
 \edm
is not. The second claim is immediately clear, since $\Psi_1$ maps $V_2$ into
$V_1$.  For the first, Theorem~\ref{action-on-invariants} implies that
$\Psi_k$ maps $I_{m}$ into a multiple of $I_{m+k}-I_{m-k}$, which 
is a linear combination of elements of $V_1$. The same applies to the
image of all differences $I_{m+k}-I_{m-k}$.
\begin{NB}
The example $G=\SL(2,\C)$ is treated in detail in Section 3 of the paper
\cite{Kostant&M01} and the authors obtain similar formulas.
\end{NB}
%
%
\section{A separation of variables theorem for $\SL(2,\C)$}
\label{s-sep-var-sl2}\noindent
%
\subsection{Harmonic cofree actions}
%

We recall that an action of a reductive group $G$ on an irreducible affine
variety $M$ is called \emph{cofree} if there exists a $G$-invariant subspace
$H$ of $\C[M]$ such that the multiplication map
 \begin{equation}
\label{separate-vars}
 H\ox \C[M]^G\lra \C[M],\quad h\ox f\lmapsto h\cdot f 
 \end{equation}
is an isomorphism of vector spaces.  Let $M//G$ be the {\em algebraic
quotient} of $M$ by $G$ (the affine variety such that $\C[M//G] \cong
\C[M]^G$), and let $\pi: M \mapsto M//G$ be the canonical projection
(see \cite{Kraft}). By using the solution to the Serre conjecture
concerning algebraic vector bundles on $\C^n$, Richardson
(\cite{Richardson81}) was able to establish a general algebraic
criterion for an action to be cofree.
\begin{thm}[Richardson]
Let $G$ be an algebraic group with reductive identity component and
$M$ a smooth irreducible affine $G$-variety. Then this action is cofree
whenever the following two conditions are satisfied:
\begin{enumerate}
\item[(i)]
 the algebra of invariants $\C[M]^G$ is a polynomial ring;
\item[(ii)]
 the fiber $\pi^{-1}(x)$ has dimension $\dim M - \dim M/\!/G$
for all $x\in M/\!/G$.
\end{enumerate}
Let $G$ be a simply connected semisimple algebraic group, $T$ a maximal
torus in $G$ and $W$ the Weyl group of $G$ relative to $T$. Then in
particular the following group actions are cofree:
\begin{enumerate}
\item[(a)] 
the conjugation action of $G$ on itself;
\item[(b)]
 the action of $W$ on $T$;
\item[(c)]
 the $K$-action on the symmetric space $G/K$, where $K$ is the fixed
point set of some involution $\theta$ of $G$.
\end{enumerate}
\end{thm}
\noindent
However, Richardson's proof gives no explicit realization of the space
$H$.

Classical results by Kostant and Kostant-Rallis (\cite{Kostant63},
\cite{Kostant&R71}) state (among others) that the isotropy
representation $\p$ of a symmetric space $G/K$ is always
cofree. Furthermore, in the factorization (\ref{separate-vars}) in
this case, the $K$-invariant subspace $H$ may always be realized as
the intersection of the kernels of a finite number of $K$-invariant
differential operators with constant coefficients, thus generalizing
the notion of harmonic polynomials for $\SO(n)$. This justifies the
following definition.
\begin{dfn}
A cofree action of a reductive algebraic group $K$ on an irreducible
affine variety $M$ will be called \emph{harmonic} if there exist
$K$-invariant differential operators $D_1,\,\ldots,D_n$ on $M$ such
that the linear space
\bdm
H\ =\ \bigcap_{i=1}^n \ker D_i
\edm
realizes the isomorphism (\ref{separate-vars}).
\end{dfn}

\begin{exa}
We start with an easy example of a harmonic Weyl group action.
\begin{thm}
The action of the Weyl group $W=\S_2$ on the maximal torus
$T\cong\C^*$ of $G=\SL(2,\C)$ is harmonic.
\end{thm}
\begin{proof}
The ring of regular functions of $T$ is isomorphic to $\C[e^z,e^{-z}]$
and the non trivial element of $\S_2$ acts hereon as the inversion
$e^{nz}\mapsto e^{-nz}$. Thus the invariant ring is exactly the
polynomial ring generated by $e^z+e^{-z}$, and one easily shows that
$\del_z$ and $\C[e^z,e^{-z}]$ together generate the ring of algebraic
differential operators on $T$. The operator
\bdm
D\ =\ (e^z-e^{-z})\del_z + (e^z+e^{-z}) \del^2_{z}
\edm
is obviously $W$-invariant and an easy calculation shows that
its kernel $H$ consists exactly of the functions $1$ and $e^z-e^{-z}$.
Since on the other hand the affine ring $\C[T]$
splits into the isotypic components of the trivial and the
signum representation, one gets
\bdm
\C[T]\ =\ 1\cdot \C[T]^W+ (e^z-e^{-z})\cdot \C[T]^W\ =\ H\ox\C[T]^W
\edm
and the action is therefore harmonic.
\end{proof}
\noindent
Notice that $(e^z-e^{-z})\del_z$ is just the $W$-invariant vector field
induced by the $W$-equivariant mapping $T\ra \h,\, h\mapsto h-h^{-1}$.
It should be possible to extend this example to wide classes of
Weyl group actions.
\end{exa}
\subsection{Harmonicity of the $\SL(2,\C)$ conjugation action}
%
The remainder of this section is devoted to the proof that the
conjugation action of $\SL(2,\C)$ on itself is harmonic. The strategy
is to guess a good candidate for the space $H$ of harmonics (this is
the easy part) and to explicitely construct a conjugation invariant
differential operator with kernel $H$.

\medskip\noindent
Under the simultaneous left and right action of $G$, the affine ring
of $\SL(2,\C)$ decomposes by Frobenius reciprocity into
 \begin{equation}
\label{peter-weyl-isom}
 \C[\SL(2,\C)]\ \cong\ \bigoplus_{d\geq 0}\ S^d V\ox (S^d V)^*\,. 
 \end{equation}
The decomposition of $\C[\SL(2,\C)]$ under the conjugation action of
$G$ then amounts to decomposing each occurring tensor product into
$G$-irreducibles. By the Clebsch-Gordon formula, we know that $S^d
V\ox (S^d V)^*=S^{2d}V\oplus S^{2(d-1)}V\oplus\cdots\oplus S^{0}V$,
where the trivial representation corresponds to the trace over $S^d
V$. Thus, we obtain
 \bdm
  \C[\SL(2,\C)]\ \cong\ \bigoplus_{d\geq 0}\ \C[\SL(2,\C)]^{\Ad N}(d)\ox
 S^d V  \,,
 \edm
where the first factor denotes the matrix functions invariant under
the lower diagonal unipotent matrices $N$ and of weight $d$. In
\cite[Satz 2.2]{Agricola01} it is shown that the infinite-dimensional
multiplicity spaces of the appearing irreducible $\SL(2,\C)$-modules
are irreducible and pairwise inequivalent modules for the canonical
action of the algebra of conju\-ga\-tion invariant differential
operators; this applies in particular to the ring of invariants itself
($d=0$). The space $ \C[\SL(2,\C)]^{\Ad N}(d)$ is spanned by the
products of the $d$-th power of the matrix coefficient $g_{12}$ with
any invariant,
 \bdm
 \C[\SL(2,\C)]^{\Ad N}(d)\ =\ (g_{12})^d\cdot \C[\SL(2,\C)]^{\Ad \SL(2,\C)}\,.
 \edm
Thus we may choose for $H$ the sum of all $\SL(2,\C)$-representations
with highest weight vector $(g_{12})^d$ for $d = 0, 1, 2, \ldots$. The
problem then is to construct an invariant differential operator $D$
with
 \bdm
 \C[\SL(2,\C)]^{\Ad N}\cap\ker D\ =\ \C[g_{12}]\,.
 \edm
\begin{NB}
The space $H$ coincides with the pull back under the map $g \mapsto
\Psi(g)$ of the harmonics on the Lie algebra $\g$, defined as in
\cite{Kostant63}. In the preprint \cite[6.2]{Kostant&M01}, it is shown
that this choice of a subspace $H$ of harmonics works for
$G=\SL(n,\C)$ for all $n$. Denote by $H(\g)$ the harmonics on $\g$,
and consider the conjugation invariant map
$\Phi_1(g)=g-\tr(g)\cdot\mathbf{1}/n$ already encountered in
Theorem~\ref{conj-inv-VF-sln}. Then Kostant and Michor prove the
isomorphism
 \bdm
 \C[G]\ \cong \ \C[G]^{\Ad G}\ox \Phi_1^* H(\g)\,.
 \edm
\end{NB}
\vspace{1ex}

We now turn to the construction of the differential operator
$D$. Besides the vector field $\Psi =\Psi_1$ from Theorem
\ref{conj-inv-VF-sl2}, the Casimir operator $\Delta$ that generates
the center of $U(\g)$ will also play a crucial role. We normalize
$\Delta$ such that
 \bdm
 \Delta\big|_{S^d V\ox (S^d V)^*}\ =\ d(d+2)\cdot\mathrm{id}\,.
 \edm
We need to keep track of the behavior of matrix functions under the
transition from the left and right to the conjugation action. For this
we will use the explicit isomorphism of
equation~(\ref{peter-weyl-isom}). Let $u,v$ be a basis of $V$ and
$x,y$ the dual basis of $V^*$. We choose the monomials $u^kv^{d-k},\,
k=0,\,\ldots, d,$ as a basis for $S^d V$ and realize the isomorphism
$S^d(V^*)\cong(S^d V)^*$ by
 \bdm
 \vphi_1 \cdots \vphi_d\lmapsto \Big[v_1 \cdots v_d \lmapsto
 \frac{1}{d!}\sum_{\sigma\in S_d}\vphi_{\sigma(1)}(v_1) \cdots
 \vphi_{\sigma(d)}(v_d) \Big]\,.
 \edm
Then $\binom{d}{k}x^ky^{d-k}$ is the basis vector in $(S^d V)^*$ dual
to $u^kv^{d-k}$. Since confusions cannot occur, we shall henceforth
omit the tensor product sign for elements of $S^d V\ox (S^d
V)^*$. With this choice of dual bases, the trace over $S^d V$ is the
total contraction and may be written
 \bdm
 \tr|_{\S^d V}\ =\ (ux+vy)^d\ =\ \sum_{k=0}^d \binom{d}{k}x^ky^{d-k}\cdot
 u^kv^{d-k}\,.
 \edm
Elements of $\SL(2,\C)$ will be parameterized as 
$g=\left[\nms\ba{ll}\alpha & \beta \\ \gamma & \delta \ea\nms\right]$. They 
act on $V$ and $V^*$ by
 \bdm
 g\cdot\left[\nms\ba{c}u\\v\ea\nms\right]\ =\ \left[\nms\ba{c}\alpha u+\beta v
\\ \gamma u + \delta v\ea\nms\right],\quad
g^{-t}\left[\nms\ba{c}x\\y\ea\nms\right]\ =\ \left[\nms\ba{c}\delta x-\gamma y
\\ -\beta x + \alpha y\ea\nms\right]\,.
 \edm
For illustration, we check that $(ux+vy)^d$ is $G$-invariant, thus
reproving its identification with the trace,
 \begin{eqnarray*}
 g\cdot (ux+vy)^d &=& \left[(\alpha u+\beta v)(\delta x-\gamma y)+
 (\gamma u + \delta v)(-\beta x + \alpha y)\right]^d 
 \\
 &=&
 \left[(\alpha\delta-\beta\gamma)(ux+vy)\right]^d = (ux+vy)^d
 \end{eqnarray*}
and compute the function on $G$ corresponding to the tensor
$(uy)^d\in S^d V\ox S^d V^*$:
 \begin{eqnarray*}
 f(g) = y^d (g\cdot u^d) &=& y^d\big((\alpha u + \beta v)^d\big) =
 y^d (\beta^d v^d + d \alpha\beta^{d-1}uv^{d-1}+\ldots)
\\
 &=&
 \beta^d\cdot 1 +  \alpha\beta^{d-1}\cdot 0+\ldots+0 = \beta^d\,.
 \end{eqnarray*}
This is just the highest weight function $(g_{12})^d$.
\begin{thm}\label{harmonicity-of-sl2-action}
The kernel of the conjugation invariant differential operator
 \bdm
 D\ =\ -\tr(g)^3 \Delta + \tr(g)\Psi^2+(\tr(g)^2+4)\Psi
 \edm
on $G=\SL(2,\C)$, intersected with $\C[\SL(2,\C)]^{\Ad N}$, consists
of the linear hull of all the functions $(g_{12})^n,\, n\in\N$. Hence
the conjugation action of $\SL(2,\C)$ is harmonic.
\end{thm}
\begin{proof}
The proof consists of a tedious computation; we only give an outline
here. As in Theorem~\ref{action-on-invariants}, one shows that the
vector field $\Psi$ acts on the matrix function $J_{m,n}(g) = J_m\cdot
g_{12}^n= (\alpha+\delta)^m \beta^n$ by
 \bdm
 \Psi(J_{m,n})\ =\ (n+m)\,J_{m+1,n}-4m\,J_{m-1,n}\,.
 \edm
Thus we obtain for the square of its action
 \bdm
 \Psi^2(J_{m,n})\ =\ (n+m)(n+m+1)J_{m+2,n}-4(2m^2+mn+n)J_{m,n}+16m(m-1)
 J_{m-2,n}\,.
 \edm
For $m=0$ this means in particular
 \bdm
 D(J_{0,n})\ =\ -\tr(g)^3n(n+2)J_{0,n}+\tr(g)\left[n(n+1)J_{2,n}-4nJ_{0,n}
 \right] + (\tr(g)^2+4)nJ_{1,n}\ =\ 0\,,
 \edm
as needed. The proof that $DJ_{m,n}\neq 0$ for $m\neq 0$ requires more
work.  The problem is to determine the function on $G$ corresponding
to the tensor $(uy)^m(ux+vy)^n \in S^{m+n}V\ox S^{m+n}V^*$, in order
to deduce the eigenvalue of $\Delta$ on $J_{m,n}$. In fact, a full
formula can only be proved for $n=0$. In this situation, one first
shows on the maximal torus $T$ of $G$ the validity of the formula
\bdm 
  \tr(g)\big|_{S^m V} \ =\ \sum_{k=0}^{[m/2]} (-1)^k \binom{m-k}{k}
     \tr(g)^{m-2k}\,.  
\edm 
Then, using $\Delta \, \tr(g)\big|_{S^m V}\ =\ m(m+2) \tr(g)\big|_{S^m
V}$, a lengthy induction proof yields
 \bdm
 \Delta \, \tr(g)^m\ = \ m(m+2)\tr(g)^m-4m(m-1)\tr(g)^{m-2}\,.
 \edm
The explicit calculation may be found in
\cite[p.\,54-55]{Agricola00}. Since the action of $\Psi$ and $\Psi^2$
was determined before, one gets for the action of $D$
 \bdm
 D \,  \tr(g)^m\ = \ -4m(m+1)\tr(g)^{m+1}+16m(m-2)\tr(g)^{m-1}\,.
 \edm
The right hand side vanishes exactly for $m=0$, as it should. For the
general case, we show that the matrix function $f_{m,n}(g)$
corresponding to $(uy)^n(ux+vy)^m$ has the form
\bdm
f_{m,n}(g)\ =\ \beta^n\left[\tr(g)^m+ \frac{m(1-m)}{n+m}\tr(g)^{m-2}+ 
R\right]\,,
\edm
where the remainder $R$ is a sum of $\tr(g)$ to the powers
$m-4,m-6,\ldots$.  The main point here is in fact the precise value of
the coefficient of $\tr(g)^{m-2}$, since the general form of this
Ansatz is obviously correct.  For $n=0$, we recover for the second
coefficient the old result $1-m=-\binom{m-1}{1}$. For the computation,
we may restrict $f_{m,n}(g)$ to the Borel subgroup $B$ of all group
elements with $\gamma=0$. Then one has for $b\in B$
 \bea[*]
 f_{m,n}(b)& = &
 (\alpha u+\beta v)^n y^n \left[ (\alpha u+\beta v)x+\alpha^{-1} vy\right]^m
 \\
 & = & \left(\beta^nv^n+n\alpha\beta^{n-1}uv^{n-1}+ \cdots \right)y^n 
 \\
 & & \cdot \left[
  \frac{v^my^m}{\alpha^m}+ m\frac{v^{m-1}y^{m-1}}{\alpha^{m-1}} 
  (\alpha u+\beta v)x 
   + \frac{m(m-1)}{2}\frac{v^{m-2}y^{m-2}}{\alpha^{m-2}}
   (\alpha u+\beta v)^2x^2 + \cdots 
   \right]\,.
 \eea[*]
We sort this product by increasing powers of $\alpha$, starting with
$\alpha^{-m}$. The product of the first summands in every factor
yields the only contribution to $\alpha^{-m}$. The two mixed products
of the first summand in one factor and the second summand in the other
factor both yield contributions to $\alpha^{-(m-2)}$ and
$\alpha^{-(m-1)}$. However, the contribution to $\alpha^{-(m-1)}$ is
zero, because $v^{n+m}xy^{n+m-1}=0$, these two basis elements are not
dual to each other. Similarly, the product of $\beta^mv^m$ in the
first factor with the third summand in the second factor gives no
contribution to $\alpha^{-(m-2)}$, because the basis vectors do not
match. To summarize, one gets the expansion
 \bdm
 f_{m,n}(b)\ =\ \beta^n\left[ \frac{v^{n+m}y^{n+m}}{\alpha^m}+
      (m+nm) \frac{uv^{n+m-1}xy^{n+m-1}}{\alpha^{m-2}} + \cdots
    \right] .
 \edm
The vector $xy^{n+m-1}$ is dual to $uv^{n+m-1}$ up to a correction
factor of $n+m$, so we finally get
 \bdm
 f_{m,n}(b)\ =\ \beta^n\left[ \frac{1}{\alpha^m}+\frac{m(n+1)}{n+m}
 \frac{1}{\alpha^{m-2}} + \cdots \right] .
 \edm
By its nature, $f_{m,n}(b)$ has to be a product of $\beta^n$ times an
invariant. Thus the expression inside the brackets is a linear
combination of powers of $(\alpha+\alpha^{-1})$.  Since
$(\alpha+\alpha^{-1})^m=\alpha^{-m}+m\alpha^{-(m-2)} + \cdots$, there
exists a rearrangement of terms such that
\bdm
 f_{m,n}(b)\ =\ \beta^n\left[(\alpha+\alpha^{-1})^m + \bigg(\frac{m(n+1)}{n+m}
 -m\bigg)(\alpha+\alpha^{-1})^{m-2} + \cdots \right]\,.
\edm
This is the desired expression, on which we can now study the action
of the Casimir operator. The function $f_{m,n}$ is an eigenfunction of
$\Delta$ with eigenvalue $(n+m)(n+m+2)$, hence
 \bdm
 \Delta \beta^n\tr(g)^m\ =\ (n+m)(n+m+2)f_{m,n}- \frac{m(1-m)}{n+m}\Delta
 \beta^n\tr(g)^{m-2}- \Delta R\,,
 \edm
and again $\Delta\beta^n\tr(g)^{m-2} =
(n+m-2)(n+m)\beta^n\tr(g)^{m-2}+$ lower order terms. Sorting by powers
of $\tr(g)$, we get
 \bdm
 \Delta \beta^n\tr(g)^m\ =\ (n+m)(n+m+2)\beta^n\tr(g)^m
 + 4m(1-m)\beta^n\tr(g)^{m-2}+\tilde{R}\,.
 \edm
We conjecture that $\tilde{R}=0$, but we do not need this here. Notice
that the second coefficient does \emph{not} depend on $n$. Going back
to the definition of the operator $D$, we sort the result again by
decreasing powers of $\tr(g)$ to obtain
 \bdm
 D \beta^n\tr(g)^m\ =\ -4m(n+m+1) \beta^n\tr(g)^{m+1} +\cdots \,.
 \edm
As a polynomial in $\tr(g)$, $D \beta^n\tr(g)^m$ vanishes if and only
if every coefficient is zero, and a look on the second factor shows
that this happens precisely for $m=0$. Thus we showed that $D
\beta^n\tr(g)^m\neq 0 $ for $m\neq 0$.
\end{proof}

%


\begin{thebibliography}{Ban-J90}

\bibitem[Agr00]{Agricola00}
I.~Agricola, \emph{Die {F}robenius-{Z}erlegung auf einer algebraischen
  ${G}$-{M}annigfaltigkeit}, Ph.D. thesis, Humboldt-Universit{\"a}t zu Berlin,
  2000.

\bibitem[Agr01]{Agricola01}
I.~Agricola, \emph{Invariante {D}ifferentialoperatoren und die
  {F}robenius-{Z}erlegung einer {$G$}-{V}ariet\"at}, J. of Lie Th.
\textbf{11} (2001), 81--109.

\bibitem[Ban-J90]{Bang-Jensen}
J.~ Bang-Jensen, \emph{The Multiplicities of Certain K-Types in Spherical
Representations}, J. Functional Analysis \textbf{91} (1990), 346--403.

\bibitem[GW97]{Goodman&W}
R.~Goodman and N.~R.~Wallach, \emph{Representations and invariants of the
  classical groups}, Encyclopedia of Mathematics and its Appl., vol.~68,
  Cambridge Univ. Press, 1998.

\bibitem[Hel78]{Helgason78} 
S.~Helgason, \emph{Differential Geometry,
Lie Groups and Symmetric Spaces}, Pure and Applied Mathematics,
vol. 80, Acad. Press, New York, 1978.

\bibitem[Hel84]{Helgason84}
S.~Helgason, \emph{Groups and geometric analysis}, Pure and Applied
  Mathematics, vol. 113, Acad. Press, Orlando, 1984.

\bibitem[Hel94]{Helgason94}
S.~Helgason, \emph{Geometric Analysis on Symmetric Spaces}, Mathematical
Surveys and Monographs, vol.~39, Amer. Math. Soc., 1994.

\bibitem[Kos63]{Kostant63}
 B.~Kostant, \emph{Lie group representations
on polynomial rings}, Amer. J. Math.  \textbf{85} (1963), 327--404.

\bibitem[KR71]{Kostant&R71}
B.~Kostant and S.~Rallis, \emph{Orbits and representations associated with
  symmetric spaces}, Amer. J. Math. \textbf{93} (1971), 753--809.

\bibitem[KM01]{Kostant&M01}
B.~Kostant and P.~Michor, \emph{The generalized Cayley map
from an algebraic group to its Lie algebra}, preprint
(arXiv:math.RT/0109066v1, 10 Sep 2001), to appear in \emph{The Orbit
  Method in Geometry and Physics} (A.A.~Kirillov Festschrift),
Progress in Mathematics, Birkh\"auser, 2003.

\bibitem[Kra85]{Kraft}
H.~Kraft, \emph{Geometrische {M}ethoden in der {I}nvariantentheorie},
  Vieweg, Braunschweig, 1985.

\bibitem[Ric81]{Richardson81}
R.~W.~Richardson, \emph{An application of the {S}erre conjecture to semisimple
  algebraic groups}, Algebra, Proc. Conf., Carbondale 1980, Lecture Notes in
  Mathematics, vol. 848, Springer (1981), 141--151.

\bibitem[Ric82]{Richardson82}
R.~W.~Richardson, \emph{Orbits, invariants and representations associated to
  involutions of reductive groups}, Invent. Math. \textbf{66} (1982), 287--312.

\bibitem[Sli-B89]{Sliman}
L.~Sliman-Bouattour, \emph{Sur les champs de vecteurs invariants sur un espace
hermitien sym\'{e}trique}, C.R.~Acad.~Sci.~Paris, \textbf{309} (1989),
439--442.

\bibitem[Sol63]{Solomon}
L.~Solomon, \emph{Invariants of finite reflection groups}, Nagoya
Math. J. \textbf{22} (1963), 57--64.

\bibitem[Ste65]{Steinberg65}
R.~Steinberg, \emph{Regular elements of semisimple algebraic groups}, Publ.
Math. I.H.E.S. \textbf{25} (1965), 49--80.

\bibitem[Vre76]{Vretare76}
L.~Vretare, \emph{Elementary spherical functions on symmetric spaces}, Math.
Scand. \textbf{39} (1976), 343--358.

\end{thebibliography}
\end{document}